\documentclass{article}

\usepackage[utf8]{inputenc}
\usepackage{fullpage}
\usepackage{amssymb}
\usepackage{amsmath}
\usepackage{graphicx}
\usepackage{wrapfig}
\usepackage{cite}

\title{\textbf{Connected and Integrated Transportation Systems}}
\author{Andreas A. Malikopoulos\\ Terri Connor Kelly and John Kelly Career Development Professor\\ University of Delaware}
\date{}

\begin{document}

\maketitle

\section{Overview}

There is a solid body of research  available that has aimed at enhancing our understanding of optimally controlling the powertrain of both conventional and electric vehicles. Many different approaches have been proposed to address the fundamental vehicle system performance challenges using both offline and online analytical algorithms.

\begin{figure}
	\centering
	\includegraphics
	[width=0.4\textwidth]
	{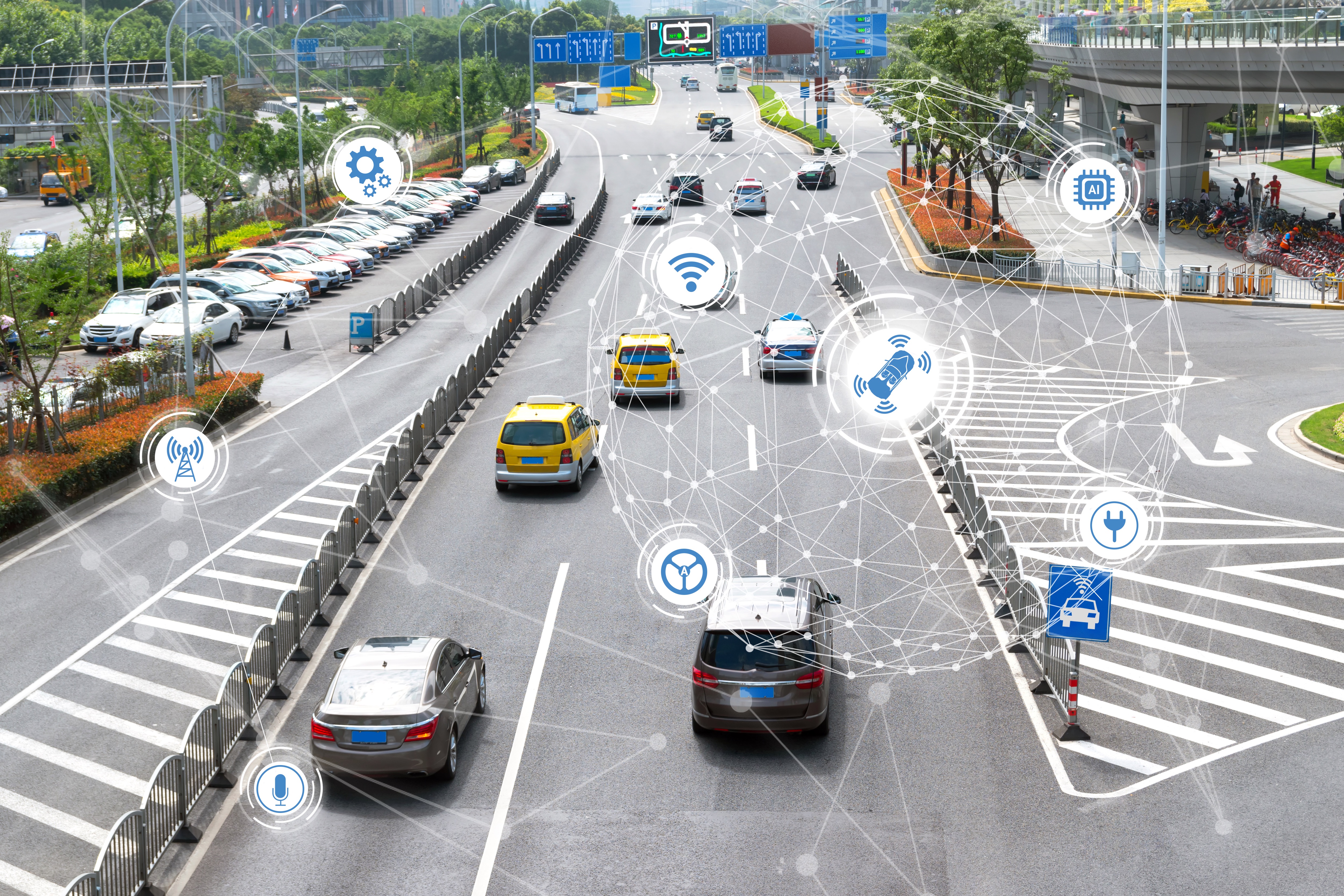}
	\caption{\small
		Connected and automated vehicles in a traffic environment.
	}
	\label{fig:1}
\end{figure}

Connected and automated vehicles (CAVs) (Fig. \ref{fig:1}) have attracted considerable attention over the last decade since they provide the most intriguing opportunity for enabling users not only to improve powertrain performance but also to better monitor transportation network conditions and make better decisions for improving safety and transportation efficiency \cite{Spieser2014, Fagnant2014}. In a transportation network with CAVs, we can consider the vehicle as part of a larger system, which can be optimized at an even larger scale. Such large-scale optimization requires the acquisition and processing of additional information from the driver and conditions outside the vehicle itself. This requires addition of new sensors and/or better utilization of information generated by existing sensors. The processing of such multi-scale information provides new opportunities for developing  significantly new approaches in order to overcome the curse of dimensionality.  It seems clear that the availability of this information has the potential to reduce traffic accidents and ease congestion by enabling vehicles to more rapidly account for changes in their mutual environment that cannot be predicted by deterministic models. Likewise, communication with traffic structures, nearby buildings, and traffic lights should allow for individual vehicle control systems to account for unpredictable changes in local infrastructure.

Recognition of the necessity for connecting vehicles to their surroundings has gained momentum. Many stakeholders intuitively see the benefits of multi-scale vehicle control systems and have started to develop business cases for their respective domains, including the automotive and insurance industries, government, and service providers. The availability of vehicle-to-vehicle (V2V) and vehicle-to-infrastructure (V2I) communication has the potential to ease congestion and improve safety by enabling vehicles to respond rapidly to changes in their mutual environment. Vehicle automation technologies can aim at developing robust vehicle control systems that can quickly respond to dynamic traffic operating conditions. With the advent of emerging information and communication technologies, we are witnessing a massive increase in the integration of our energy, transportation, and cyber networks. These advances, coupled with human factors, are giving rise to a new level of complexity \cite{Malikopoulos2016c} in transportation networks. As we move to increasingly complex emerging transportation systems, with changing landscapes enabled by connectivity and automation, future transportation networks could shift dramatically with the large-scale deployment of CAVs. On the one hand, with the generation of massive amounts of data from vehicles and infrastructure, there are opportunities to develop optimization methods to identify and realize a substantial energy reduction of the transportation network, and to optimize the large-scale system behavior using the interplay among vehicles. 

\begin{figure}
	\centering
	\includegraphics
	[width=0.6\textwidth]
	{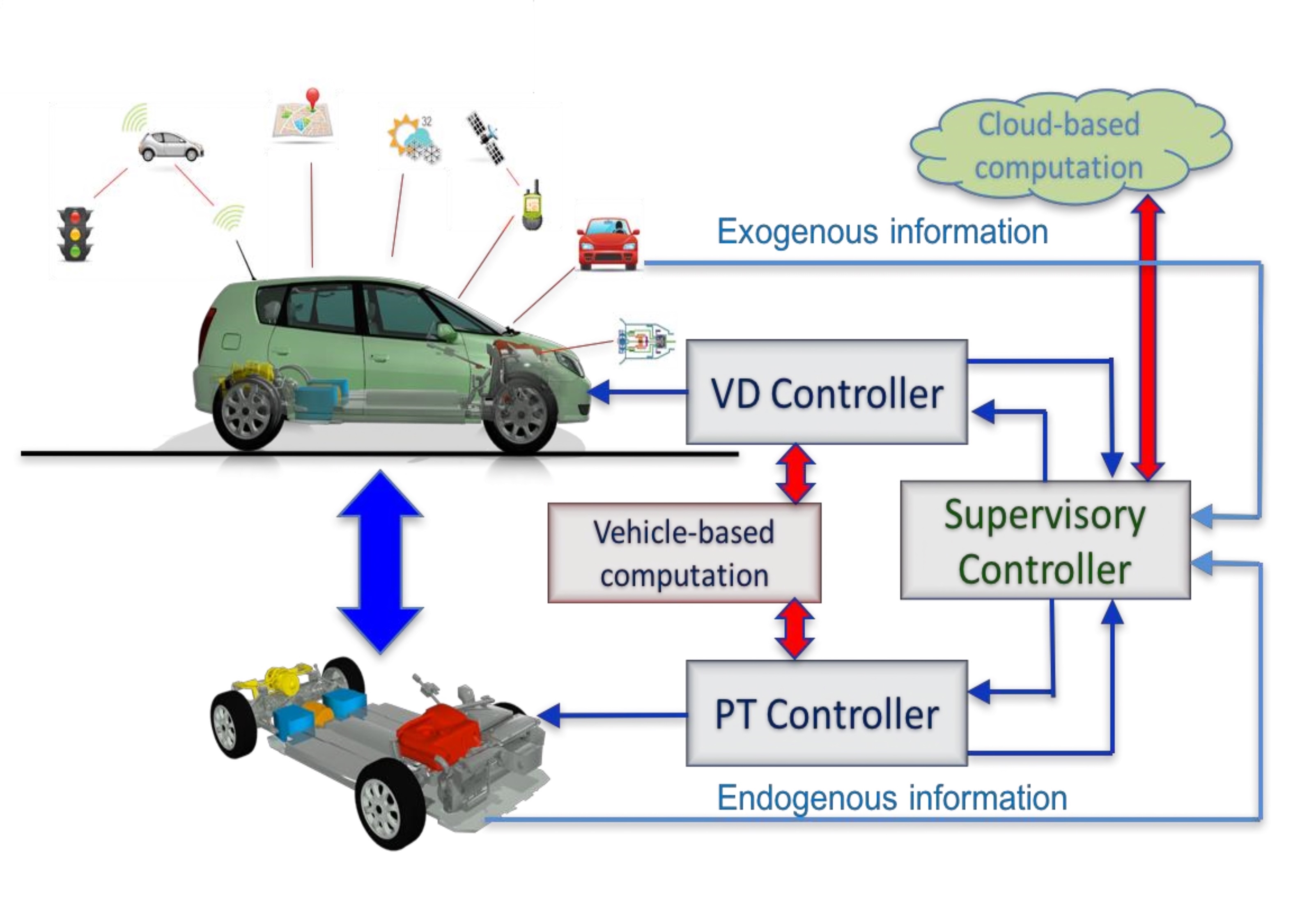}
	\caption{\small
		Control architecture for a connected and automated vehicle.
	}
	\label{fig:2}
\end{figure}

In a transportation network with CAVs there is additional information available that can be used to control and optimize jointly both vehicle-level and powertrain-level operation. For example, control technologies have been reported  aimed at maximizing the energy efficiency of a 2016 Audi A3 e-tron plug-in hybrid electric vehicle by more than 25\% without degradation in tailpipe out exhaust emission levels, and without sacrificing the vehicle’s drivability, performance, and safety \cite{mahbub2020sae-2}. These  technologies can: (I) optimize the vehicle’s speed profile aimed at minimizing (ideally, eliminating) stop-and-go driving, and (II) optimize the powertrain of the vehicle for this optimal speed profile obtained under (I). The control control architecture of such technologies (Fig. \ref{fig:2}) consist of (a) a vehicle dynamic (VD) controller, (b) a powertrain (PT) controller, and (c) a supervisory controller.
\begin{enumerate}
	\item The supervisory controller (1) oversees the VD and PT controllers, (2) communicates the endogenous and exogenous information appropriately, (3) computes the optimal routing for any desired origin-destination, (4) determines the regions where electric driving will have the major impact to derive a desired battery state-of-charge trajectory, and (5) synthesizes a description of the upcoming road segment from the exogenous information that communicates it to the VD controller. 
	\item The VD controller optimizes online the acceleration/deceleration and speed profile of the vehicle, and thus, the vehicle’s torque demand. 
	\item The PT controller computes the optimal nominal operation (“set-points”) for the engine, motor, battery, and transmission corresponding to the optimal solution of the VD controller. 
\end{enumerate}

The optimal solution of the VD controller along with the endogenous and exogenous information is communicated to the PT controller through the supervisory controller. A unique feature of the control architecture is that the supervisory controller coordinates the VD and PT controllers to ensure the optimal solution yielded by the VD controller is feasible for the PT controller, and eventually results in maximization of the vehicle’s energy efficiency. If the optimal solution of the VD controller is not feasible, then the supervisory controller enforces the VD controller to repeat the optimization for a new set of parameters. 
To this end, we focus on approaches that have been developed for the VD controller.

\subsection{Vehicle Dynamic Controller}

In a typical commute, we encounter traffic scenarios that include crossing intersections, merging at roadways and roundabouts,  cruising in congested traffic, passing through speed reduction zones, and lane-merging or passing maneuvers.
There have been two major approaches to utilizing connectivity and automation to improve transportation efficiency and safety, namely, (1) platooning and (2) traffic smoothing. 

The first approach utilizes connectivity and automation to form closely-coupled vehicular platoons to  reduce aerodynamic drag effectively, especially at high cruising speeds. The concept of forming platoons of vehicles traveling at a high speed 
was a popular system-level approach to address traffic congestion, which gained momentum in the 1980s and 1990s \cite{Shladover1991,Rajamani2000}. Such automated transportation system can alleviate congestion, reduce energy use and emissions, and improve safety while increasing throughput significantly. The Japan ITS Energy Project \cite{tsugawa2013}, the Safe Road Trains for the Environment program \cite{davila2010sartre}, and the California Partner for Advanced Transportation Technology \cite{shladover2007path} are among the mostly-reported efforts in this area.

The second approach is to smooth the traffic flow through a VD controller for optimal coordination of CAVs at different traffic scenarios, e.g, at merging roadways \cite{Rios-Torres2017,Ntousakis2016aa}, intersections \cite{Malikopoulos2017,Malikopoulos2020,Dresner2008,Lee2012a,gregoire2014priority,fayazi2018mixed}, adjacent intersections \cite{chalaki2020TCST,chalaki2020TITS}, speed reduction zones \cite{Malikopoulos2018c}, roundabouts \cite{chalaki2020experimental}, and 
corridors \cite{Zhao2018ITSC, mahbub2020decentralized,zhao2019CCTA-2,chalaki2021CSM}.  One of the very early efforts in this direction was proposed by Athans \cite{Athans1969} for safe and efficient coordination of merging maneuvers with the intention of avoiding congestion. Assuming a given merging sequence, the merging problem was formulated as a linear optimal regulator \cite{Levine1966} to control a single string of vehicles, with the aim of minimizing the speed errors that will affect the desired headway between each consecutive pair of vehicles.  In 2004, Dresner and Stone \cite{Dresner2004} proposed the use of the reservation scheme to control a single intersection of two roads. In their approach, each vehicle requests the reservation of the space-time cells to cross the intersection at a particular time interval defined from the estimated arrival time to the intersection.  Since then, numerous approaches have been proposed on coordinating CAVs to improve traffic flow \cite{Kachroo1997, Antoniotti1997, Ran1999}, and to achieve safe and efficient control of traffic through various traffic bottlenecks where potential vehicle collisions may happen \cite{Dresner2008,DeLaFortelle2010, Huang2012,Zohdy2012,Yan2009,Li2006,Zhu2015a,Wu2014,kim2014}.   Queuing theory has also been used to address this problem by modeling coordination of CAVs as a polling system with two queues and one server that determines the sequence of times assigned to the vehicles on each road \cite{Miculescu2014}.  Some of the methods presented in the literature have focused on multi-objective optimization problems  \cite{Kamal2013a, Kamal2014, Campos2014, Makarem2013, qian2015}. More recently, a study \cite{Ratti2016} indicated that transitioning from intersections with traffic lights to autonomous intersections, where vehicles can coordinate and cross the intersection without the use of traffic lights, has the potential of doubling capacity and reducing delays.
Two survey papers that report the research efforts in this area can be found in \cite{Malikopoulos2016a, guanetti2018}.

\subsection{Impact and Future Directions}

Several research efforts have focused on quantifying the impact of CAVs on vehicle miles traveled, energy, and greenhouse gas (GHG) emissions \cite{Martin2011,Martin2011a,firnkorn2015,Martin2016,Chen2016}. Some studies \cite{Martin2011,Martin2011a,Martin2016} have shown a decrease in GHG emissions  with significant implications on public transit \cite{Martin2011a}. Other research efforts have investigated the feasibility and potential environmental impacts \cite{Shaheen2014} of shared CAVs \cite{chong2013a,Ford2012,Rigole2014,Bischoff2016,Dia2017,Merlin2017,Metz2018,Fiedler2018,Lu2018}. There have been also studies focusing on cost-benefit analysis of a mobility system with  CAVs \cite{Burns2012,Fagnant2014,Bosch2017,Moorthy2017} and the impact on vehicle ownership by using surveys or comparable analysis with conventional car-sharing systems \cite{Truong2017,tussyadiah2017,Foldes2018,Hao2018,Menon2018}. There are several survey papers providing a good review in related topics, e.g., see \cite{Jorge2013,Agatz2012a,Furuhata2013a,Brandstatter2016,lavieri2017,jittrapirom2017,utriainen2018}.

\begin{figure}
	\centering
	\includegraphics
	[width=0.7\textwidth]
	{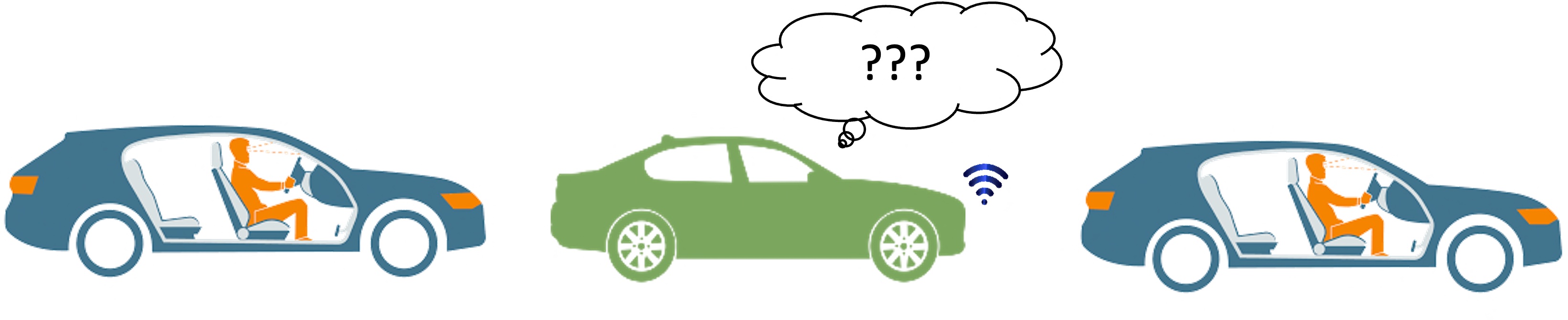}
	\caption{\small
		A mixed traffic environment consisting of human-driven and connected automated vehicles.
	}
	\label{fig:3}
\end{figure}

Recent studies \cite{Malikopoulos2018d, Zhao2018CTA}, have investigated the impact of different penetration rates of CAVs, e.g., 0\% to 100\%, in fuel consumption and travel time for two traffic scenarios: (1) merging at highways on ramp and (2) merging at roundabouts. What was observed is that as the penetration rate of CAVs is decreased, both fuel consumption and travel time deteriorate. 
It is expected that CAVs will gradually penetrate the market, interact with human-driven vehicles (HDVs) (Fig. \ref{fig:3}), and contend with V2V and V2I communication limitations, e.g., bandwidth, dropouts, errors and/or delays, as well as perception delays, lack of state information, etc. However, different levels of vehicle automation  in the transportation network can significantly alter  transportation efficiency metrics  ranging from 45\% improvement to 60\% deterioration \cite{wadud2016}.  Moreover, we anticipate that efficient  transportation and travel cost reduction might alter human travel behavior causing rebound effects, e.g., by improving efficiency, travel cost is decreased, hence willingness-to-travel is increased. The latter would increase overall vehicle miles traveled, which in turn might negate the benefits in terms of energy and travel time.

As we move to increasingly diverse mobility systems with different penetration rates of CAVs, new approaches are needed to optimize the impact on system behavior of the interplay between CAVs and HDVs at different traffic scenarios.  While several studies have shown the benefits of CAVs to reduce energy and alleviate traffic congestion in specific traffic scenarios, most of these efforts have focused on $100$\% CAV penetration rates without considering HDVs. One key question that still remains unanswered is ``how can CAVs and HDVs be coordinated safely?''

While several studies discussed in this chapter have shown the benefits of emerging mobility systems to reduce energy and alleviate traffic congestion in specific transportation scenarios, the research community should focus on developing a mobility system that can enhance accessibility, safety, and equity in transportation without causing rebound effects, while also gaining the travelers' acceptance. Future research should address this question and attempt to design a ``socially-optimal mobility system,'' i.e.,  a mobility system that (1) is efficient (in terms of energy consumption and travel time), (2) mitigates rebound effects, and (3) ensures equity in transportation.\\

\bibliographystyle{IEEEtran}
\bibliography{TAC_Ref_Andreas,TAC_Ref_Andreas2}

\end{document}